\documentclass[12pt]{amsart}
\usepackage{amsmath,latexsym,amscd,amsbsy,amssymb,amsfonts,amsthm,fleqn,leqno,
euscript, graphicx, texdraw}
\numberwithin{equation}{section}

\newcommand{\cl}{\mathrm{cl}}

\newcommand{\I}{\mathbb I}

\newtheorem{thm}{Theorem}

\newtheorem{lem}[thm]{Lemma}

\theoremstyle{definition}

\newtheorem{re}[thm]{Remark}

\theoremstyle{remark}

\begin{document}

\title[Mazurkiewicz manifolds and homogeneity]
{Mazurkiewicz manifolds and homogeneity}

\author{P. Krupski}
\address{Mathematical Institute, University of Wroc\l aw, pl.
Grunwaldzki 2/4, 50--384 Wroc\l aw, Poland}
\email{Pawel.Krupski@math.uni.wroc.pl}

\author{V. Valov}
\address{Department of Computer Science and Mathematics,
Nipissing University, 100 College Drive, P.O. Box 5002, North Bay,
ON, P1B 8L7, Canada} \email{veskov@nipissingu.ca}
\thanks{The second author was partially supported by NSERC Grant 261914-08.}

\date{\today}


\begin{abstract}
It is proved that no region of a homogeneous locally compact,
locally connected metric space can be cut by an $F_\sigma$-subset of
a ``smaller" dimension. The result applies to different finite or
infinite topological dimensions of metrizable spaces.
\end{abstract}
\maketitle\markboth{}{}

The classical Hurewicz-Menger-Tumarkin theorem in dimension theory
says that connected topological $n$-manifolds (with or without
boundary) are Cantor manifolds (i.~e., no subset of covering
dimension $\le n-2$ separates the space). The theorem was almost
immediately strengthened by Mazurkiewicz who proved that  regions
(i. e., open connected subsets) in Euclidean spaces (and, in fact,
in topological manifolds) cannot be cut  by subsets of codimension
at least two (a subset \emph{cuts} if its complement is not
continuum-wise connected~\cite{E}).   The  Hurewicz-Menger-Tumarkin
theorem has many generalizations. In particular, it is known that
regions of homogeneous locally compact metric spaces are Cantor
manifolds (including their infinite-dimensional
versions)~\cite{Kr,Kr1}. It was proved in~\cite{HT} that no weakly
infinite-dimensional subset cuts the  product
 of a countable number of nondegenerate metric continua.

In this paper, we obtain a generalization in  spirit of the
Mazurkiewicz theorem: regions in  homogeneous locally compact,
locally connected  metric spaces cannot be cut by
$F_\sigma$-subsets of  codimension at least two. Moreover, our
result holds true for a very general dimension function
$D_{\mathcal{K}}$ considered in~\cite{KKTV}  which captures the
covering dimension, cohomological dimension $\dim_G$ with respect to
any Abelian group $G$ as well as the extraordinary dimension
$\dim_L$ with respect to a given $CW$-complex $L$, and has its
counterparts in infinite dimensions including $C$-spaces and weakly
infinite-dimensional spaces.

Basic facts on Cantor manifolds  and their stronger variations with
respect to dimension $D_{\mathcal{K}}$ or to the above-mentioned
infinite dimensions has been presented in~\cite{KKTV}. We recall
some necessary terminology and results from that paper. We restrict
our considerations to metrizable spaces.

A sequence $\mathcal{K}=\{K_0,K_1,..\}$ of $CW$-complexes is called
a {\em stratum} for a dimension theory~\cite{dr} if
\begin{itemize}\item
for each space $X$ admitting a perfect map onto a metrizable space,
$K_n\in AE(X)$ implies both $K_{n+1}\in AE(X\times\I)$ and
$K_{n+j}\in AE(X)$ for all $j\geq 0$.
\end{itemize}
Here, $K_n\in AE(X)$ means that $K_n$ is an absolute extensor for
$X$. Given a stratum $\mathcal{K}$, the dimension function
$D_{\mathcal{K}}$ for a metrizable  space $X$ is defined as follows:
\begin{enumerate}
\item
$D_{\mathcal{K}}(X)=-1$ iff $X=\emptyset$;
\item $D_{\mathcal{K}}(X)\le n$ if
$K_n\in AE(X)$ for $n\ge 0$; if $D_{\mathcal{K}}(X)\le n$ and
$K_m\not\in AE(X)$  for all $m<n$, then $D_{\mathcal{K}}(X)= n$;
\item
 $D_{\mathcal{K}}(X)=\infty$ if $D_{\mathcal{K}}(X)\le n$ is not satisfied for any $n$.
\end{enumerate}

According to the countable sum theorem in extension theory, it
follows directly from the above definition that
$D_{\mathcal{K}}(X)\leq n$ implies $D_{\mathcal{K}}(A)\leq n$ for
any $F_{\sigma}$-subset $A\subset X$.

\

Henceforth,  $\mathcal{C}$ will denote one of the four classes of
metrizable spaces:
\begin{enumerate}
\item
 the class $\mathcal D_{\mathcal K}^k$ of at most
$k$-dimensional spaces with respect to dimension $D_{\mathcal{K}}$,

\item the class $\mathcal D_{\mathcal{K}}^{<\infty}$ of strongly countable
$D_{\mathcal K}$-dimensional spaces, i.e. all spaces represented as a countable union
of closed finite-dimensional subsets with respect to
$D_{\mathcal{K}}$,

\item
the class $\mathbf C$ of  $C$-spaces,

 and
\item the class $\mathcal {WID}$
of weakly infinite-dimensional spaces.
\end{enumerate}

(for definitions of a weakly (strongly) infinite-dimensional or of a $C$-space, see~\cite{E}).

\

A metrizable space $X$  is a {\em  Cantor manifold with respect to a
class $\mathcal{C}$} if $X$ cannot be separated by a closed subset
which belongs to $\mathcal{C}$.

$X$  is a {\em Mazurkiewicz manifold with respect to  $\mathcal{C}$}
if for every two closed, disjoint subsets $X_0,X_1\subset X$, both
having non-empty interiors in $X$, and every $F_\sigma$-subset
$F\subset X$ with $F\in\mathcal{C}$, there exists a continuum in
$X\setminus F$ joining $X_0$ and $X_1$.

Obviously, Mazurkiewicz manifolds with respect to  $\mathcal{C}$
are  Cantor manifolds with respect to $\mathcal{C}$.
 It was observed in~\cite{KKTV} that  if no
$F_\sigma$-subset from a class $\mathcal{C}$ cuts a
compact space $X$, then $X$ is a Mazurkiewicz manifold with respect
to $\mathcal{C}$; the converse implication holds for locally
connected compact spaces $X$.

\

We are going to use the following theorem from~\cite{KKTV}.

\begin{thm}\label{exist}
Let $Z$ be a metric compact space and $Z\notin\mathcal{C}$, where
$\mathcal{C}$ is any of the  following four classes:
$\mathcal{WID}$, $\mathbf C$, $\mathcal D_{\mathcal K}^{n-2}$,
$\mathcal D_{\mathcal K}^{<\infty}$. In the case
$\mathcal{C}=\mathcal D_{\mathcal K}^{n-2}$   we additionally assume
$D_{\mathcal{K}}(Z)=n$ and in the case $\mathcal C=\mathcal
D_{\mathcal K}^{<\infty}$ assume that $Z$ does not contain closed
subsets of arbitrary large finite dimension $D_{\mathcal K}$. Then
$Z$ contains a Mazurkiewicz manifold with respect to  $\mathcal{C}$.
\end{thm}

We also need the following version of the Effros' theorem
(see~\cite[Proposition 1.4]{Kr1}).

\begin{thm}\label{effros}
If $X$ is a homogeneous locally compact metric space  (with metric
$\rho$), then for every $a\in X$ and $\epsilon>0$ there exists
$\delta>0$ such that if $\rho(x,a)<\delta$, then there is an
$\epsilon$-homeomorphism $h:X\to X$ (i.~e., $\rho(h(y),y)<\epsilon$
for each $y$) such that $h(a)=x$.
\end{thm}

The following Lemma is a slight generalization of~\cite[Theorem 8,
p.~243]{Ku}.

\begin{lem}\label{lemma}
If $X$  is a locally compact, locally connected metric space  and
the union  $\bigcup_{i=1}^\infty F_i$  cuts a region $U$ of $X$,
where $F_i$ is a closed subset of $U$ for each $i$, then there is
$i$ such that $F_i$ cuts a region $V\subset U$.
\end{lem}
\begin{proof}
Choose two distinct points $a,b\in U$ such that
$\bigcup_{i=1}^\infty F_i$ cuts $U$ between them. Suppose no set
$F_i$ cuts any subregion of $U$. So, there is a subregion
$U_1\subset U\setminus F_1$ containing $a$ and $b$. Since $U_1$ is
completely metrizable, it is arcwise connected (by the
Mazurkiewicz-Moore-Menger theorem~\cite{Ku}), hence there is an arc
$\alpha_1\subset U_1$ from $a$ to $b$.  The local compactness and
local connectedness allows us to get a region $U_1'$ such that
$\alpha_1\subset U_1'\subset \cl(U_1')\subset U_1$ and $\cl(U_1')$
is compact. Similarly, we find an arc   $\alpha_2$ from $a$ to $b$
and regions $U_2\subset U_1'\setminus F_2$ and $U_2'$ such that
$\alpha_2\subset U_2'\subset \cl(U_2')\subset U_2$ and $\cl(U_2')$
is compact. Continuing this way, we get a decreasing sequence of
continua $\cl(U_n')$ whose intersection is a continuum in
$U\setminus \bigcup_{i=1}^\infty F_i$  containing $a$ and $b$, a
contradiction.
  \end{proof}

  \

We can now prove our main result.

\begin{thm}\label{main}
Let $X$ be a homogeneous locally  compact, locally connected metric
space. Suppose $U$ is a region in $X$ and $U\notin \mathcal C$,
where $\mathcal C$ is one of the above four classes. In case
$\mathcal{C}=\mathcal D_{\mathcal K}^{n-2}$ assume
$D_{\mathcal{K}}(U)=n$. Then $U$ is a Mazurkiewicz manifold with
respect to  $\mathcal{C}$.
\end{thm}

\begin{proof}
Notice that $U$ is second countable. It follows (by the countable
sum theorem for spaces in class  $\mathcal{C}$) that $U$ contains
compact sets of arbitrary small diameters which do not belong to
$\mathcal{C}$. Suppose  $U$ is not a Mazurkiewicz manifold  with
respect to  $\mathcal{C}$  and let  an $F_\sigma$-subset
$\bigcup_{i=1}^\infty F_i$  of $U$ cut $U$ with each $F_i\in\mathcal
C $ being closed in $U$. It follows by Lemma~\ref{lemma} that there
is $j$ such that $F_j$ cuts a region $V\subset U$. Thus, $F_j$ also
separates $V$ since $V$ is locally connected~\cite[Theorem 1, p.
238]{Ku}. Without loss of generality one can assume that $F_j$ is
nowhere dense. Let $V\setminus F_j=V_1\cup V_2$, where $V_1$ and
$V_2$ are nonempty, open and disjoint.  Fix a point $a\in F_j\cap
V$.

Assume first that in the case $\mathcal C=\mathcal D_{\mathcal
K}^{<\infty}$ the region $U$ does not contain closed subsets of
arbitrary large finite dimension $D_{\mathcal K}$. Then, by
Theorem~\ref{exist}, there are  arbitrary small compact Mazurkiewicz
manifolds  with respect to  $\mathcal{C}$ in $U$. By the
homogeneity, there is a compact Mazurkiewicz manifold (with respect
to $\mathcal{C}$)  $M\subset V$ containing $a$.  $M$ being  a Cantor
manifold with respect to $\mathcal{C}$, it is not in $\mathcal{C}$,
so $M$ is not contained in $F_j$. Suppose $M$ intersects $V_1$. Then
the Effros' Theorem~\ref{effros} allows us to  push $M$ toward $V_2$
by a small homeomorphism  so that it meet both sets $V_1$ and $V_2$.
This means that the displaced $M$ is separated by $F_j$, a
contradiction.

The case when  $\mathcal C=\mathcal D_{\mathcal K}^{<\infty}$ and
$U$  contains closed subsets of arbitrary large finite dimension
$D_{\mathcal K}$ can be handled in a similar way. Indeed,  let
$D_{\mathcal K}(F_j)=m$ and observe  that, by  the
$\sigma$-compactness of $U$ and by the countable sum theorem for
dimension $D_{\mathcal K}$, we can assume that these closed subsets
are compact and of arbitrary small diameters.    By
Theorem~\ref{exist}   they contain arbitrary small compact
Mazurkiewicz manifolds with respect to  corresponding finite
dimensions $D_{\mathcal K}$. Thus there are arbitrary small compact
Cantor manifolds with respect to $\mathcal D_{\mathcal K}^{k-2}$ in
$U$  for some $k> m+2$. By the homogeneity, there is such  a  Cantor
manifold $M\subset U$ containing $a$. Now, using the Effros'
Theorem~\ref{effros},  we get a contradiction as in the previous
paragraph.
\end{proof}

\begin{re}
If $\mathcal C=\mathcal D_{\mathcal K}^{<\infty}$, then the
hypothesis $U\notin \mathcal C$ in Theorem~\ref{main} can be
equivalently replaced  by $D_{\mathcal K}(U)=\infty$. More
precisely, we have the following proposition (cf.~\cite[Proposition
4.3]{KKTV}).

\emph{If $X$ is a homogeneous locally compact metric space and $U$
is a second countable open subset of $X$, then  $D_{\mathcal
K}(U)=\infty$ if and only if $U\notin \mathcal D_{\mathcal
K}^{<\infty}$.}

Indeed, suppose $D_{\mathcal K}(U)=\infty$ but
$U=\bigcup_{i=1}^\infty F_i$, where each $F_i$ is a closed  subset
of $U$ of finite dimension $D_{\mathcal K}$. Since each closed
subset of $U$ can be represented as a countable union of compact
subsets, we can assume that each $F_i$ is compact. Then the Baire
theorem and  the homogeneity easily imply that $U$ is contained in a
union of countably many homeomorphic copies of some $F_{i_0}$. So,
as an $F_\sigma$ subset of the union, its dimension $D_{\mathcal
K}(U)$ is finite, a contradiction. The converse implication is
obvious.
\end{re}

\end{document}